\documentclass[a4paper]{article}
\usepackage{amsmath}
\usepackage{amsfonts}
\usepackage{amsthm}
\usepackage{amssymb}
\usepackage[USenglish]{babel}
\usepackage[utf8]{inputenc} %para traducir acentos
\usepackage[T1]{fontenc}
\usepackage{tensor}
\usepackage{times}
\usepackage{cite}

\usepackage{parskip} % to avoid indentation of all paragraph
\allowdisplaybreaks

\newtheorem{theorem}{Theorem}[section]
\newtheorem{proposition}{Proposition}[section]
\newtheorem{question}{Question}

\theoremstyle{remark}
\newtheorem{remark}{Remark}[section]

\begin{document}

\title{A note on Kirchhoff's theorem for almost complex spheres I}

\author{L\'azaro O. Rodr\'iguez D\'iaz \thanks{IMPA, Rio de Janeiro, Brazil; email:lazarord@impa.br}}

\date{}

\maketitle

\begin{abstract}
By a theorem of Kirchhoff if the six sphere admits an almost complex structure then the seven sphere is parallelizable, more crucial, he exhibited an explicit global frame constructed out of the given almost complex structure. This result implicitly equips the seven sphere with a definite $H$-space multiplication. We propose to address the existence problem of complex structures on the six sphere studying the associated parallelism-multiplications on the seven sphere. We ask to what extent the integrability condition of the almost complex structure amounts to the constancy of the structure functions of the global frame defining the parallelism, i.e, if this parallelism comes from a Lie group structure. At a more fundamental level we inquire if the integrability condition of the almost complex structure entails the homotopy associativity of the induced multiplication. A positive answer to these questions would rule out the six sphere of being a complex manifold since the seven sphere is not a Lie group, not even a homotopy associative $H$-space.
\end{abstract}

\section{Introduction} 

It is well known the only spheres that are parallelizable are $S^1$, $S^{3}$ and $S^{7}$, in all cases it is possible to exhibit a trivialization of the tangent bundle using the complex numbers, quaternions and octonions respectively. However, among the spheres, $S^1$ and $S^{3}$ are the only Lie groups, again it is the structure of the complex numbers and quaternions which induces the Lie group structures. A natural candidate group structure for $S^{7}$ using the octonions fails due to the non-associativity of the octonions, see Remark \ref{rm-parallelism}.

 Another interesting geometric structure to consider on spheres are almost complex structures, in this case $S^{0}, S^{2}$ and $S^{6}$ are the only ones that admit an almost complex structure \cite{Serre-Borel}. There are others proofs of this fact, one of them combines a result of Kirchhoff \cite{Kirchhoff} that if $S^{n}$ carries an almost complex structure then $S^{n+1}$ is parallelizable and the list of parallelizable manifolds, i.e., $S^1$, $S^{3}$ and $S^{7}$, \cite{BottMilnor58, Kervaire58}. Summed up: \emph{$S^{n}$ admits an almost complex structure if and only if $S^{n+1}$ is parallelizable.}

The octonions induces an almost complex on the sphere $S^{6}$, nevertheless it is not integrable \cite{Eckmann1951,Ehresmann51}. We can prove this is due to the non-associativity of the octonions, see Remark \ref{rm-integrability}. It is unknown whether $S^{6}$ admits an integrable almost complex structure.    

It is natural to ask if the probable non-existence of integrable almost complex structure on $S^{6}$ is related to $S^{7}$ not being a Lie group, or more generally, if it is linked to the non-existence of homotopy associative multiplications on $S^{7}$. After all, the natural candidates to integrable almost complex structure on $S^6$ and Lie group structure on $S^7$ fail for the same reason, non-associativity of the octonions. 

As mentioned above it was proved by Kirchhoff \cite{Kirchhoff} that if the sphere $S^{n}$ admits an almost complex structure $J$, then $S^{n+1}$ is parallelizable, more crucial, he exhibited the explicit absolute parallelism, i.e., an explicit global frame $\{X_{1}, \cdots, X_{n}\}$ constructed using $J$. This result implicitly equips the seven sphere with a definite $H$-space multiplication.

We propose to address the existence problem of complex structures on the six sphere studying the associated parallelism-multiplications on the seven sphere. We ask to what extent the integrability condition of the almost complex structure amounts to the constancy of the structure functions of the global frame defining the parallelism, i.e, if this parallelism comes from a Lie group structure. At a more fundamental level we inquire if the integrability condition of the almost complex structure entails the homotopy associativity of the induced multiplication. A positive answer to these questions would rule out the six sphere of being a complex manifold since the seven sphere is not a Lie group, not even a homotopy associative $H$-space.

\section{When a parallelism comes from a Lie group structure?}
 
To each smooth global frame $\{X_{1}, \cdots, X_{n}\}$ on a manifold $M$ there is associated a  flat (zero curvature) connection $\Gamma$.
 The covariant derivative of $\Gamma$ is given by $\nabla_{Z}\left(\sum f^{i}X_{i}\right)=\sum Z(f^{i})X_{i}$, for any smooth vector field $Z$ on $M$. 

The structure equations of $\Gamma$ in the frame $\{X_{1}, \cdots, X_{n}\}$ are:
\begin{gather}\label{eq-MaurerCartan}
d\theta^{i}=\tfrac{1}{2}\tensor{T}{^i_j_k}\theta^{j}\wedge \theta^{k} \quad \text{and} \quad \omega^{i}_{j}=0, 
\end{gather}
\noindent
where $\{\theta^{1},\cdots, \theta^{n}\}$ is the coframe dual of $\{X_{1}, \cdots, X_{n}\}$ and $w^{i}_{j}$ are the connection forms. The torsion tensor of $\Gamma$ is given by
\begin{equation}\label{eq-structurefunctions}
T(X_{j},X_{k})=\sum_{i=1}^{n}\tensor{T}{^{i}_{jk}}X_{i}=-[X_{j},X_{k}].
\end{equation}

Furthermore the torsion tensor of $\Gamma$ is parallel if and only if the \emph{structure functions} $\tensor{T}{^i_j_k}$ are constant. More generally a tensor field on $M$ is parallel with respect to $\Gamma$ if and only if has constant components with respect to the field of frames $\{X_{1}, \cdots, X_{n}\}$.

The equation (\ref{eq-structurefunctions}) resembles the way we define the structure constants of a Lie algebra, i.e., the first equation in (\ref{eq-MaurerCartan}) looks like the Maurer-Cartan equation;  this is not coincidence since Lie groups are always parallelizable. Moreover the following converse states when a parallelism on $M$ comes from a Lie group structure on $M$:
\vspace{0.1in}
\begin{theorem}[\cite{Hicks59}, Theorem $5$]\label{thm-Hicks}
Let $M$ be a simply connected manifold on which is defined a complete linear connection with zero curvature and torsion invariant under parallel
translation. Then M admits a Lie group structure such that left translations induce the original connection.
\end{theorem}

The above theorem was proved by Chern \cite[Section 5, page 128]{Chern53} formulated in terms of an $\{e\}$-structure on $M$. Generalization of this result to not necessarily simply connected manifolds has been proved many times in the literature, e.g., Wolf \cite[Proposition 2.5]{Wolf72} within the context of absolute parallelism. Certainly all these is rooted in Cartan's local equivalence method, see Sternberg \cite[Theorem 2.4, Chapter V]{Sternberg64}.

In general a linear connection on a compact manifold is not necessarily complete. However as the geodesics of the associated connection $\Gamma$ consist of the integral curves of the vectors fields $\{X_{1},\cdots,X_{n}\}$, compactness of $M$ implies $\Gamma$ is complete.

In a local coordinate patch $U$ of $M$ with coordinates $(x_{1},\cdots,x_{n})$, we can write the vector fields $X_{j}=\sum X_{j}^{l}\frac{\partial}{\partial x^{l}}$ and the dual forms $\theta^{j}=\sum \theta_{l}^{j}dx^{l}$, (where $\left(\theta^{j}_{l}\right)$ is the inverse matrix of $\left(X_{j}^{l}\right)$ ) in terms of the local basis of vectors fields and differential forms. One way to compute the structure functions is using the formula:
\begin{equation}\label{eq-structure_functions}
\tensor{T}{^i_j_k}=\sum_{r,s=1}^{n}X^{r}_{j}X^{s}_{k}\left(\frac{\partial \theta^{i}_{s}}{\partial x^{r}}-\frac{\partial \theta^{i}_{r}}{\partial x^{s}}\right).
\end{equation}

\subsection{Classical parallelism of the seven sphere}

Let briefly discuss the parallelism of $S^{7}$ induced by the octonions in the light of Hicks's theorem, Theorem \ref{thm-Hicks}. The next proposition is a generalization of how orthogonal multiplications are useful to define vector fields on spheres.
\vspace{0.1in}
\begin{proposition}[\cite{brickell1970}, Proposition $7.3.1$]\label{thm-Brickell}
Suppose we have a map $\nu:  \mathbb{R}^{k+1}\times\mathbb{R}^{n+1}\rightarrow\mathbb{R}^{n+1}$  linear in the first factor and continuous in the second factor satisfying:
\begin{itemize}
\item[i)] $\nu(v,z)=0$ implies $z=0$ or $v=0$,
\item[ii)] there exists $e\in\mathbb{R}^{k+1}$ such that $\nu(e,z)=z$ for all $z\in\mathbb{R}^{n+1}$,
\end{itemize}
then $S^{n}$ admits $k$ independent vector fields.
\end{proposition}

The proof can be found in \cite{brickell1970}, there was stated for bilinear maps $\nu$, however the proof only uses linearity in the first factor. As a corollary of the above proposition the multiplications in the complex, quaternions and octonions numbers induce a parallelism in $S^{1}, S^{3}$ and $S^{7}$ respectively. 

Fix the canonical basis of the octonions $\mathbb{O}$  given by the identity $1$ and seven imaginary octonions $e_{i}$, $i=1,\cdots, 7$ satisfying the multiplication rule: $e_{i}e_{j}=-\delta_{ij}+a_{ijk}e_{k}$, where the structure constants $a_{ijk}$ are totally antisymmetric in the three indices. Using Proposition \ref{thm-Brickell} we construct seven linearly independent vector fields $X_{i}$ on the sphere $S^{7}\subset \mathbb{O}$ of octonions of norm one as follow: $X_{i}(x)=e_{i}x$ for $x\in S^{7}$, $i=1,\cdots, 7$. 

Lets compute the structure functions of this global frame. Note the multiplication in this particular case is linear in both factors, therefore the Lie brackets $[X_{i},X_{j}]$ can be computed by the commutator of the corresponding linear maps.
\begin{align}\label{eq_Torsion_octonions}
 [X_{i},X_{j}](x)&=e_{i}(e_{j}x)-e_{j}(e_{i}x)\\\nonumber
&=2a_{ijk}e_{k}x-2[e_{i},e_{j},x]\\
&=2\left(a_{ijk}- \langle[e_{i},e_{j},x], e_{k}x \rangle \right) X_{k}(x),\nonumber
\end{align} 
where $[a,b,c]:=(ab)c-a(bc)$ is the associator,  $\langle a, b\rangle:=\tfrac{1}{2}\left(a\bar{b}+b\bar{a}\right)$ is the scalar product and the conjugation is defined by $\bar{1}=1$, $\bar{e_{i}}=-e_{i}$ and $\overline{ab}=\bar{b}\bar{a}$.
\vspace{0.1in}
\begin{remark}\label{rm-parallelism}
\emph{Thus the non-associativity of the octonions causes the non-constancy of the structure functions of the classical parallelism of $S^{7}$}. Compare to Remark \ref{rm-integrability}. The non-commutativity of the algebra causes the non-vanishing of the torsion.
\end{remark}
\vspace{0.1in}
\begin{remark}\label{rm-alternativity-torsion}
Note the structure functions coincide with the structure constants of the algebra at the north and south pole, i.e., at $1$ and $-1$. Compare to Remark \ref{rm-Steenrod_approximation_theorem}.
\end{remark}
\vspace{0.1in}
\begin{remark}
We used the alternativity of the octonionic product to prove the second equality in (\ref{eq_Torsion_octonions}). Compare to Remark \ref{rm-alternativity-Nijenhuis}.
\end{remark}

\section{An integrability condition in disuse}

An important theorem of Newlander and Nirenberg \cite{NewlanderNirenberg57} states an almost complex structure $J$ on a manifold $M$ is integrable if and only if the Nijenhuis tensor $N(X,Y)=[JX,JY]-[X,Y]-J[X,JY]-J[JX,Y]$ vanishes identically, where $X$ and $Y$ are vector fields on $M$. 

For any torsion free connection $\nabla$ on $M$ we have $N(X,Y)=\left(\nabla_{JX}J-J\nabla_{X}J\right)Y-\left(\nabla_{JY}J-J\nabla_{Y}J\right)X$ . 
In local coordinates $(x_1, \cdots, x_{2n})$ in $M$, the components $N_{jk}^{i}$ of $N$ can be expressed in terms of the components $J^{i}_{j}$ of $J$ and its partial derivatives:
\begin{align*}
N^{i}_{jk}=&\sum_{l=1}^{2n}\left( J^{l}_{j}\left(\frac{\partial J^{i}_{k}}{\partial x^{l}}-\frac{\partial J^{i}_{l}}{\partial x^{k}}\right)-J^{l}_{k}\left(\frac{\partial J^{i}_{j}}{\partial x^{l}}-\frac{\partial J^{i}_{l}}{\partial x^{j}}\right)\right).
\end{align*}
There are others less used equivalent formulations of the integrability condition of an almost complex structure. Lets recall one which seems suitable for the problem at hand. See for example Calabi-Spencer \cite{calabi_spencer1951}, Guggenheimer \cite[Theorem 13]{Guggenheimer52}, Hodge \cite[page 105]{Hodge52}, a good account appears also in Calabi \cite[section 3]{Calabi58}.

Define the folllowing action of $J$ on the algebra of differential forms: $J(f):=f$ for $f\in C^{\infty}(M)$ and $\left(Jw\right)\left(X_{1},\cdots,X_{p}\right):=w\left(JX_{1},\cdots,JX_{p}\right)$, for $w$ a $p$-form on $M$. In other words it is the action of the dual map of $J$ on $1$-forms extended to higher order forms by the above formula.
\vspace{0.1in}
\begin{proposition}\label{thm-calabi}
An almost complex structure $J$ on $M$ is integrable if and only if
\begin{equation}\label{eq-calabi}
\left(dJdJ-JdJd\right)f=0,
\end{equation}
is satisfied for all functions $f$ of class $C^{2}$.
\end{proposition}

Moreover, Hodge proved that Proposition \ref{thm-calabi} is also true if $\left(dJdJ-JdJd\right)w=0$ is satisfied for $p$-forms $w$ of degree $1\leq p \leq \operatorname{dim}M-2 $, \cite[page 105]{Hodge52}. In local coordinates the above proposition said:
\vspace{0.1in}
\begin{proposition}\label{thm-calabi-spencer}
An almost complex structure $J$ on $M$ is integrable if and only if the tensor 
\begin{equation}\label{eq-calabi-spencer}
\tau^{i}_{jk}=\sum_{p,q=1}^{2n}\left(\delta^{p}_{j}\delta^{q}_{k}- J^{p}_{j}J^{q}_{k}\right)\left(\frac{\partial J^{i}_{q}}{\partial x^{p}}-\frac{\partial J^{i}_{p}}{\partial x^{q}}\right),
\end{equation}
vanishes identically.
\end{proposition}
It is clear that $\tau^{i}_{jk}=-J^{r}_{j}N^{i}_{rk}$,  however, Proposition \ref{thm-calabi} shows the vanishing of the tensor $\tau^{i}_{jk}$ can be interpreted as the integrability condition on its own. 
\vspace{0.1in}
\begin{remark}
Bearing in mind $J^{-1}=-J$, note the analogy between (\ref{eq-calabi-spencer}) and the formula to compute the structure functions (\ref{eq-structure_functions}).
\end{remark}

\section{Kirchhoff's theorem and the octonions}

Before going into Kirchhoff's theorem let briefly recall how the multiplication in the octonions $\mathbb{O}$ induces an almost complex structure on $S^{6}$. Kirchhoff construction is modeled on this, in fact its proof reverses this process, he reconstructs the `multiplication' of $\mathbb{R}^{8}$ from the almost complex structure, see Remark \ref{rm-Brickell} and Section \ref{sec-H_structures}.
\subsection{The almost complex structure induced by the octonions}
Let $\mathop{Im}\mathbb{O}\subset \mathbb{O}$ denotes the hyperplane of imaginary octonions orthogonal to $1\in\mathbb{O}$ and let $S^{6}\subset \mathop{Im}\mathbb{O}$ be the sphere of imaginary octonions of norm one. Right multiplication by $y\in S^{6}$ induces an orthogonal linear transformation $R_{y}:\mathbb{O}\rightarrow \mathbb{O}$ that satisfies $\left(R_{y}\right)^{2}=-1$. Moreover, $R_{y}$ preserves the plane spanned by $1$ and $y$ ($1\rightarrow y$, $y\rightarrow -1$), therefore preserves its orthogonal six dimensional plane, which can be identified with $T_{y}S^{6}\subset \mathbb{O}$. It follows $R_{y}$ induces an almost complex structure on $S^{6}$.  Now we are going to show  the Nijenhuis tensor corresponding to this almost complex structure can be written in terms of the associator of $\mathbb{O}$. 

The Nijenhuis tensor can be computed by:
\begin{align*}
N(X,Y)&=d(JY)(JX)-d(JX)(JY)-dY(X)+dX(Y)\\
&-J\left(d(JY)(X)-dX(JY)\right)-J\left(dY(JX)-d(JX)(Y)\right),
\end{align*}
to see this, note we are in euclidean space, then we can compute the Lie brackets of two vector fields  $X:S^{6}\rightarrow \mathbb{R}^{7}$, $Y:S^{6}\rightarrow \mathbb{R}^{7}$ by $[X,Y]=dY(X)-dX(Y)$, where $dX$ and $dY$ denote the differential of $X$ and $Y$ respectively as maps.

By definition $J_{a}Y_{a}=Y_{a} \cdot a $ where $ a\in S^{6}$ and $Y$ is a vector field on $S^{6}$, differentiating we get:
\begin{align*}
d(JY)(JX)=&J\left(dY(JX)\right) + Y\cdot JX,\\ J\left(d(JY)(X)\right)=&(-1)dY(X)+J(Y\cdot X).
\end{align*}
Then
\begin{equation*}
N(X,Y)=Y\cdot JX- X\cdot JY-J(Y\cdot X)+J(X\cdot Y).
\end{equation*}
For $b, c \in T_{a}S^{6}$ we get:
\begin{equation}\label{eq_Nijenhuis_associator}
N_{a}(b,c)=c\cdot (b\cdot a)- b\cdot (c\cdot a)-(c\cdot b)\cdot a+(b\cdot c)\cdot a=2[a,b,c].
\end{equation}
\begin{remark}\label{rm-integrability}
 \emph{Therefore the non-associativity of the octonions is responsible for the non-integrability of this almost complex structure}. Compare to Remark \ref{rm-parallelism}.
\end{remark}
\vspace{0.1in}
\begin{remark}\label{rm-alternativity-Nijenhuis}
To establish the last equality in (\ref{eq_Nijenhuis_associator}) we used the algebra of octonions is alternative. Compare to Remark \ref{rm-alternativity-torsion}.
\end{remark}

\subsection{Kirchhoff's theorem}

Now we state Kirchhoff's theorem and its proof \cite{Kirchhoff}, \cite[Theorem V]{Kirchhoff53}, see also Kobayashi-Nomizu \cite[Chapter IX, Example 2.6]{Kobayashi_NomizuII}.
\vspace{0.1in}
\begin{theorem}[\cite{Kirchhoff}, Theorem $4$]\label{thm-Kirchhoff_theorem}
If the sphere $S^{n}$ admits an almost complex structure, then $S^{n+1}$ is parallelizable.
\begin{proof}
We need to exhibit a field $\sigma$ of linear frames on $S^{n+1}$.
Let $J$ be an almost complex structure on $S^{n}$. In the vector space $\mathbb{R}^{n+2}$ fix a subspace $\mathbb{R}^{n+1}$ and a unit vector $e:=e_{n+2}\in \mathbb{R}^{n+2}$ perpendicular to $\mathbb{R}^{n+1}$. Denote by $S^{n}$ and  $S^{n+1}$  the unit spheres in $\mathbb{R}^{n+1}$ and $\mathbb{R}^{n+2}$ respectively.

Given $y\in S^{n}$ denote by $V_{y}$ the $n$-dimensional vector subspace of $\mathbb{R}^{n+2}$ parallel to the tangent space $T_{y}(S^{n})$ in $\mathbb{R}^{n+2}$ and $J_{y}$ the linear endomorphism of $V_{y}$ corresponding to the linear endomorphism of $T_{y}(S^{n})$ given by $J$.
Define a linear transformation $\widetilde{J}_{y}:\mathbb{R}^{n+2}\rightarrow \mathbb{R}^{n+2}$ by $\widetilde{J}_{y}(e)=y$, $\widetilde{J}_{y}(y)=-e$ and $\widetilde{J}_{y}(z)=J_{y}(z)$ for $z\in V_{y}$.  It follows from  $J^{2}=-\mathop{Id}$ that ${\widetilde{J}_{y}}^{2}=\mathop{-Id}$.

 Let $x \in \mathbb{R}^{n+2}$, then it can be written uniquely as follows:
\begin{gather}\label{eq-decomposition}
x=\alpha e + \beta y, \quad \alpha, \beta \in \mathbb{R}, \quad \beta\geq 0,  \quad \text{and} \quad y \in S^{n}.
\end{gather} 
Define the linear transformation:
\begin{gather}\label{eq-linearframe}
\widetilde{\sigma}_{x}: \mathbb{R}^{n+2}\rightarrow\mathbb{R}^{n+2}, \quad
\widetilde{\sigma}_{x}:=\alpha \mathop{Id}+\beta \widetilde{J}_{y},
\end{gather}
where $\mathop{Id}$ denotes the identity transformation of $\mathbb{R}^{n+2}$.  As ${\widetilde{J}_{y}}^{2}=\mathop{-Id}$ we get that $\widetilde{\sigma}_{x}$ is an isomorphism. Note also that $\widetilde{\sigma}_{x}(e)=x$.
If the transformations $\widetilde{\sigma}_{x}$ for $x\in S^{n+1}$, are restricted to $\mathbb{R}^{n+1}$,
\begin{gather*}
\sigma_{x}:=\left.\widetilde{\sigma}_{x}\right|_{\mathbb{R}^{n+1}}, \quad  x\in S^{n+1},
\end{gather*}
 we get the desire linear frame $\sigma_{x}=:\mathbb{R}^{n+1}\rightarrow T_{x}(S^{n+1})$. In fact, $\mathbb{R}^{n+1}$ is spanned by $y$ and $V_{y}$ and that both $\sigma_{x}(y)$ and $\sigma_{x}(z)$, $z\in V_{y}$ are perpendicular to $x$ in $\mathbb{R}^{n+2}$, then can be considered as elements of $T_{x}(S^{n+1})$.
\end{proof}
\end{theorem}
\begin{remark}
Note that Kirchhoff's theorem does not assume any additional condition on  the almost complex structure $J$.
\end{remark}
\vspace{0.1in}
\begin{remark}\label{rm-almost_Hermitian}
If we assume $J$ is an almost hermitian structure, i.e., $J$ is  compatible with some Riemannian metric $g$ on $S^{n}$ (this is always possible), the theorem above was rewritten by Steenrod \cite[Theorem 41.19]{Steenrod51}. He noted that in this case the constructed global frame $\widetilde{\sigma}$ is in fact an orthogonal frame, that is, $\widetilde{\sigma}_{x}\left(\widetilde{\sigma}_{x}\right)^{T}=Id$, for $x\in S^{n+1} $. More generally, we have $\widetilde{\sigma}_{x}\left(\widetilde{\sigma}_{x}\right)^{T}=\|x\|^{2}Id$, when $x\in \mathbb{R}^{n+2}$. This will be used in Section \ref{sec-H_structures}. In what follows $\|x\|$ will denote the Euclidean norm of $x$.
\end{remark}
\vspace{0.1in}
\begin{remark}
The vector fields $\{X_{i}(x):=\sigma_{x}(e_{i})\}_{i=1,\cdots,n+1}$ defining the parallelism in Theorem \ref{thm-Kirchhoff_theorem} can be written explicitly as:
\begin{equation*}
X_{i}(x)=x_{n+2}e_{i}-x_{i}e_{n+2}+\beta(x)J_{y}\left(e_{i}-\langle y, e_{i}\rangle y\right),
\end{equation*} 
where $\{e_{i}\}_{i=1,\cdots,n+2}$ is the canonical basis of  $\mathbb{R}^{n+2}$.
\end{remark}
\vspace{0.1in}
\begin{remark}\label{rm-Steenrod_approximation_theorem}
Note the linear frame $\sigma$ is smooth at all points of $S^{n+1}$ except at $e$ and $-e$, where it is only continuous. Steenrod's approximation theorem \cite[Theorem 6.7]{Steenrod51} tell us every continuous section of a smooth fibre bundle can be approximated arbitrarily closely by a smooth section, even more, preserving the original section on regions where it is already differentiable. Therefore we will assume $\sigma$ is a smooth linear frame that coincides on $S^{n+1}\setminus\{e,-e\}$ with the constructed in the above theorem.   
\end{remark}
\vspace{0.1in}
\begin{remark}\label{rm-inverse}
If $x\in S^{n+1}$ then $\alpha$ and $\beta$ in (\ref{eq-decomposition}) satifies $\alpha^{2}+\beta^{2}=1$, consequently it is easy to get the inverse of (\ref{eq-linearframe}), i.e., $\widetilde{\sigma_{x}}^{-1}=\alpha\mathop{Id}-\beta\widetilde{J}_{y}$.
\end{remark}
\vspace{0.1in}
\begin{remark}\label{rm-Brickell}
We can use Proposition \ref{thm-Brickell} to conclude $S^{n+1}$ is parallelizable in Kirchhoff's theorem, for define $\nu(v,z):=\widetilde{\sigma}_{z}v$, it was proved along the proof $\nu$ satisfies the hypothesis of the Proposition \ref{thm-Brickell}.
\end{remark}

\subsection{The first approach}

The ideas discussed so far lead to formulate the following  question. Assume $S^{6}$ admits an integrable almost complex structure $J$, in particular, $J$ is an almost complex structure. By Kirchhoff's Theorem \ref{thm-Kirchhoff_theorem} the sphere $S^{7}$ is parallelizable, moreover, we have the explicit form of the global frame defining the parallelism in terms of $J$.
\vspace{0.1in}
\begin{question}\label{question1}
To what extent the integrability condition of the almost complex structure $J$ on the six sphere amounts to the constancy of the structure functions of the global frame defining the parallelism of the seven sphere?
\end{question}

Even if in general the structure functions of the frame are not constant, we can compute the complete set of invariants functions for this $\{e\}$-structure on $S^{7}$ in terms of the almost complex structure $J$. Note the torsion tensor is never zero. 

Keeping the notation used in the proof of Theorem \ref{thm-Kirchhoff_theorem}, we have seven linearly independent vector fields  $\{\sigma\left(e_{1}\right),\cdots, \sigma\left(e_{7}\right)\}$ on $S^{7}$, where $\{e_{1},\cdots, e_{7}\}$ denotes the canonical basis of $\mathbb{R}^{7}$. As we are in euclidean space, we can consider each vector field of this frame as a map $\sigma(e_{i}): S^{7}\rightarrow \mathbb{R}^{8}$, $\sigma(e_{i})(x):=\sigma_{x}(e_{i})$.
To prove the structure functions are constant is equivalent to show that $\sigma^{-1}\left([\sigma(e_{i}),\sigma(e_{j})]\right)$ is constant for each $i$ and $j$, that is,
\begin{align*}
d\left(\sigma^{-1}\left([\sigma(e_{i}),\sigma(e_{j})]\right)\right)=d\left(\sigma^{-1}\left(d\left(\sigma(e_{j})\right)\left(\sigma(e_{i})\right)-d\left(\sigma(e_{i})\right)\left(\sigma(e_{j})\right)\right)\right)=0.
\end{align*}
Observe the integrability of $J$ in the form  $\left(dJdJ-JdJd\right)w=0$ of equation (\ref{eq-calabi}) should be helpful here.

\section{H-space structures induced by almost complex structures}\label{sec-H_structures}

We can rephrase Kirchhoff's theorem as follows: if $S^{n}$ admits an almost complex structure $J$ then $S^{n+1}$ is a an $H$-space (i.e., a space which admits a continuous multiplication with a two-sided identity element). This is trivial at glance because it is well known a parallelizable sphere is an $H$-space \cite{Adams60}. The point is that the induced multiplication on $S^{n+1}$ is written explicitly in terms of $J$. 

In what follows we keep the notation of Theorem \ref{thm-Kirchhoff_theorem}. Let $J$ be an almost complex structure on $S^{6}$, and let $\widetilde{\sigma}$ be the global frame given by Kirchhoff's theorem. Define the map:
\begin{gather*}
m: S^{7}\times S^{7}\longrightarrow S^{7}, \quad m(x,y):=\widetilde{\sigma}_{x}\left(y\right)/\|\widetilde{\sigma}_{x}\left(y\right)\|.
\end{gather*}
It follows from Kirchhoff's theorem this map is a well defined multiplication on $S^{7}$ with $e$ as two-sided identity. Thus every almost complex structure $J$ on $S^{6}$ defines an $H$-space structure on $S^{7}$. 

Moreover, recall that we can always assume the given almost complex structure $J$ is compatible with some Riemannian metric $g$ on $S^{7}$. Then using Remark \ref{rm-almost_Hermitian} we can define a multiplication in $\mathbb{R}^{8}$:
\begin{gather*}
\hat{m}: \mathbb{R}^{8}\times \mathbb{R}^{8}\longrightarrow \mathbb{R}^{8}, \quad \hat{m}(x,y):=\widetilde{\sigma}_{x}\left(y\right),
\end{gather*}
which satisfies the norm product rule $\|\hat{m}(x,y)\|^{2}=\|x\|^{2}\|y\|^{2}$, and has a two-sided identity $e$. The multiplication $\hat{m}$ restricts to the $H$-space multiplication $m$ on $S^{7}$.

By a celebrated theorem of Adams \cite{Adams60} the only spheres that admit an $H$-space structure are $S^{0}, S^{1}$, $S^{3}$ and $S^{7}$. The next theorem of Wallace shows $S^{7}$ does not admit an associative multiplication.
\vspace{0.1in}
\begin{theorem}\cite[Corollary 2, page 48]{Wallace53}
	If a compact manifold $M$ admits a continuous, associative multiplication with identity then it is a topological group.	
\end{theorem}

Moreover, combining the above theorem with von Neumann's solution of Hilbert’s Fifth problem for compact groups \cite{Neumann33} would imply $M$ has a Lie group structure.

The stronger result proved in this direction is James's theorem:
\vspace{0.1in}
\begin{theorem}\cite[Theorem 1.4]{James57_II} 
	There exists no homotopy-associative multiplication on $S^{n}$ unless $n=1$ or $3$.
\end{theorem}

As we have already seen in (\ref{eq_Nijenhuis_associator}), the non-associativity of the octonions causes the non-integrability of the almost complex structure induced on $S^{6}$ by the octonions. We would like to relate the probable non-existence of complex structure on $S^{6}$ with the lack of associative or more generally homotopy associative multiplications on $S^{7}$.
 
\subsection{The second approach}
\begin{question}\label{question2}
Does the integrability condition of the almost complex structure $J$ implies the associativity or more generally the homotopy associativity of the multiplications $m$ and $\hat{m}$? 
\end{question}

Working with the multiplication $\hat{m}$ instead of $m$ has the advantage we can make use of the additive structure of $\mathbb{R}^{8}$.

As observed in Remark \ref{rm-alternativity-Nijenhuis} we employed the alternativity of the octonions in establishing the connection between the Nijenhuis tensor and the associator. It would be reasonable to expect the answer to Question \ref{question2} would require something similar to this. In this regard, it is interesting to mention that Norman proved the following:
\vspace{0.1in}
\begin{theorem}\cite[Corollary 9.3]{Norman63}
Any multiplication on a sphere satisfies the Moufang law $\left(x\cdot y\right)\left(z\cdot x\right)=\left(x\left(y\cdot z\right)\right)x$ up to homotopy.
\end{theorem}
\vspace{0.1in}
\begin{remark}
In case $J$ be the almost complex structure of $S^{6}$ induced by the octonions, Question \ref{question1} is related to Question \ref{question2}. This happens because $J$ comes a priori from an ambient multiplication of $\mathbb{R}^{8}$, in other words, if we interpret $J$ as a map $\widetilde{J}:S^{6}\longrightarrow \mathop{SO}(8,\mathbb{R})$, $x \longmapsto \widetilde{J}_{x}$  it extends as a linear map between $\mathbb{R}^{8}$ and $\mathop{SO}(8,\mathbb{R})$.     
\end{remark}
\vspace{0.1in}
\begin{remark}
However, in the general case of an almost complex structure on $S^{6}$ Question 1 and Question 2 stand at different levels and are not correlated. Question \ref{question1} asks if it possible to integrate the parallelism, i.e., if there exists a multiplication whose differential essentially induces the global frame. Question \ref{question2} asks directly if the multiplication induced by the global frame is associative (homotopy associative).
\end{remark}

\bibliographystyle{abbrv}
\bibliography{bibliografia}

\end{document}